\documentclass[12pt]{amsart}
\usepackage[latin1]{inputenc}
\usepackage{latexsym}
\usepackage{amsfonts}
\usepackage{amsmath,amssymb}
\usepackage{amssymb}
\usepackage{amscd}
\usepackage[all]{xy}

\theoremstyle{plain}
\newtheorem{theorem}{Theorem}[section]
\newtheorem{lemma}[theorem]{Lemma}
\newtheorem{proposition}[theorem]{Proposition}
\newtheorem{corollary}[theorem]{Corollary}

\theoremstyle{definition}

\theoremstyle{remark}
\newtheorem*{remark}{Remark}

\begin{document}
\title[Characterizations of trivial maps]{Characterizations of trivial maps in 3-dimensional real Milnor fibers}

\author{T. Souza, M. Hohlenwerger, D. De Mattos, R. Araújo dos Santos}

\address{Taciana O. Souza: Universidade de S\~ao Paulo - ICMC - Av. Trabalhador S\~ao-Carlense, 400 - Centro, Postal Box 668, 13560-970, S\~ao Carlos - S\~ao Paulo - Brazil \\ e-mail:tacioli@icmc.usp.br}

\address{Maria A. P. Hohlenwerger: Universidade de S\~ao Paulo - ICMC - Av. Trabalhador S\~ao-Carlense, 400 - Centro, Postal Box 668, 13560-970, S\~ao Carlos - S\~ao Paulo - Brazil \\ e-mail:amelia@icmc.usp.br}

\address{Denise De Mattos: Universidade de S\~ao Paulo - ICMC - Av. Trabalhador S\~ao-Carlense, 400 - Centro, Postal Box 668, 13560-970, S\~ao Carlos - S\~ao Paulo - Brazil \\ e-mail:deniseml@icmc.usp.br}

\address{Raimundo Araújo dos Santos: Universidade de S\~ao Paulo - ICMC - Av. Trabalhador S\~ao-Carlense, 400 - Centro, Postal Box 668, 13560-970, S\~ao Carlos - S\~ao Paulo - Brazil \\ e-mail:rnonato@icmc.usp.br}




\begin{abstract} In this paper we extend the characterization of trivial map-germs for the real Milnor fibrations started by Church and Lamotke in \cite{CL}. Our main result cover all cases on the three dimensional real Milnor fibers.

\end{abstract}

\maketitle


\section{Introduction}

In \cite{Mi} John W. Milnor showed that given a representative of a
holomorphic function germ $\psi: U\subset \mathbb{C}^{n+1} \to
\mathbb{C}$ with, $\psi(0)=0$, there exists a small enough real
number $\epsilon_{0}>0$, such that for all $0<\epsilon \leq \epsilon_{0}$

\begin{equation*}\label{eq:milnormap2}
\frac{\psi}{\|\psi \|} : S_{\epsilon}^{2n+1}\setminus K_{\epsilon}
\to S^{1}
\end{equation*}

is a smooth projection of a locally trivial fiber bundle, where
$K_{\epsilon}=\psi^{-1}(0)\cap S_{\epsilon}^{2n+1}$ is called the
link of the singularity at origin. It is well known that in the
complex setting and $n\geq 1$ the link $K_{\epsilon}$ is never empty.
Denote $F_{\theta}=\psi^{-1}(e^{i\theta})$ the fiber of the
fibration above, where $e^{i\theta}\in S^{1}$. Using tools of Morse theory,
Milnor proved that $F_{\theta}$ is a paralelizable $2n-$dimensional real manifold
and has the homotopy type of a finite
$CW-$complex of dimension $n$. Moreover, its topological closure is given by
$\overline{F_{\theta}}=F_{\theta}\cup K_{\epsilon}$
and the topological space $K_{\epsilon}$ is
$(n-2)-$connected, i.e.
$\pi_{i}(K_{\epsilon})=0$ for all $i=0,\ldots, n-2.$
It means that, for $n=2$ the link is connected and for $n\geq 3$ it
is simple connected.

\vspace{0.2cm}

In the case where $0\in\mathbb{C}^{n+1}$ is an isolated singular
point of $\psi $, Milnor gave more details about the topology of the fiber and the link.
He proved that the fiber $F_{\theta}$
has the homotopy type of a wedge or bouquet of $n-$dimensional spheres
$\displaystyle{\underbrace{S^{n}\vee \cdots \vee S^{n}}_{\mu}}$ ,
where the number of spheres $\mu$ is given by the topological degree,
$\deg_{0}\left(\displaystyle{\frac{ \nabla \psi }{\| \nabla \psi\|}}\right),$ of the mapping

$$\displaystyle{\frac{ \nabla \psi }{\| \nabla \psi \|}:S_{\epsilon}^{2n+1}\to S^{2n+1}}.$$
From this follow easily that the Euler-Poincaré
characteristic of the fiber is given by $\displaystyle{\chi(F_{\theta})=1+(-1)^{n}\mu }.$
Milnor also proved that for $\epsilon>0$ small enough, the manifold
$(f^{-1}(0)-\{0\})\cap B_{\epsilon}^{2n+2}$
intersect transversally all spheres $S_{\epsilon}^{2n+1}$ and so by
Tranversality Theorem the space $K_{\epsilon}$ is a
$(2n-1)-$dimensional smooth manifold.
Furthermore, for each $\theta$ the space $\overline{F_{\theta}}=F_{\theta}\cup K_{\epsilon}$
is a compact manifold whose boundary is the link $K_{\epsilon}$.

\vspace{0.2cm}

In the real settings, Milnor considered a real polynomial mapping
$f: \mathbb{R}^{n}\to \mathbb{R}^{p}$, $f(0)=0$, $\displaystyle{n\geq p \geq 2},$ and assumed
that in some open neighbourhood $U$ of the origin $0\in
\mathbb{R}^{n}$ we have $\Sigma(f)\cap U\subseteq \{0\}$,
where $\Sigma(f)=\{x\in U: \text{rank $(Jf)(x)$ fails to be
maximal} \}$, i.e., $0$ is an isolated singular point of $f$. Next result follow from \cite{Mi}, Theorem 11.2, page 97:

\begin{theorem}

There exists $\epsilon_{0} >0$ small enough such that, for all
$0<\epsilon \leq \epsilon_{0},$ there exists $\eta$, $0<\eta \ll
\epsilon $, such that

\begin{equation}
\displaystyle{f_{|}:f^{-1}(S_{\eta}^{p-1})\cap
B_{\epsilon}^{n}\to S_{\eta}^{p-1}}
\end{equation}

is a smooth projection of a locally trivial fiber bundle, where
$B_{\epsilon}^{n}$ stand for the $n-$dimensional closed ball
centered at origin.

\end{theorem}

It follows from $(1)$ above that the fiber is a smooth $(n-p)-$dimensional compact manifold
with boundary given by $\displaystyle{(f_{|})^{-1}(y)\cap S^{n-1}_{\epsilon}}$. Moreover, this fibration induces the following smooth fiber bundle

\begin{equation}
\displaystyle{f_{|}:f^{-1}(S_{\eta}^{p-1})\cap
\stackrel{\circ}{B_{\epsilon}^{n}}\to S_{\eta}^{p-1}}
\end{equation}

where $\stackrel{\circ}{B_{\epsilon}^{n}}$ denote the interior of the closed ball,
i.e. the $n-$dimensional open ball of radius $\epsilon $. In this case, the fiber is an open manifold.

\vspace{0.2cm}

If the link $K_{\epsilon}$ is not empty, then it is isotopic to the boundary of the fiber of fibration $(1)$. In fact, since $0\in \mathbb{R}^{n}$ is an isolated singular point of $f$,
then we have that for each $\epsilon>0$ small enough, there exists $\eta$,
$0<\eta \ll \epsilon $, such that the manifold $(f^{-1}(0)-\{0\})\cap B_{\epsilon}^{n}$
intersect transversally all spheres $S_{\epsilon}^{n-1}.$ Since the base space is contractible, the mapping

\begin{equation}
\displaystyle{f_{|}:f^{-1}(B_{\eta}^{p})\cap
S_{\epsilon}^{n-1}\to B_{\eta}^{p}}
\end{equation}

\vspace{0.2cm}

\noindent is a smooth trivial fiber bundle. Furthermore, as an application of the Ehresmann Theorem for manifolds with boundary it is not difficult to show that the mapping

\begin{equation}
\displaystyle{f_{|}:f^{-1}(B_{\eta}^{p}- \{0\})\cap
B_{\epsilon}^{n}\to B_{\eta}^{p}- \{0\}}
\end{equation}

is a smooth projection of a locally trivial fiber bundle, where
$B_{\eta}^{p}- \{0\}$ denotes the punched $p-$dimensional closed
ball in $\mathbb{R}^{p}$ of radius $\eta.$ Hence, since the inclusion $S_{\eta}^{p-1}\subset B_{\eta}^{p}-\{0\}$ is a
homotopy equivalence then, up to homotopy, the topological
information contained in the fibrations $(1)$ and $(4)$ are the
same. In addition, using the fact that a compact manifold $M$ with no empty
boundary $\partial{M}$ is homotopy equivalent to its interior
$\stackrel{\circ}{M}=M-\partial M$ then, up to homotopy, the information contained in the
fibers of fibrations $(1)$ and $(2)$ are the same. For these reasons, we can
concentrate the topological study on the fibration $(1).$

\begin{remark}

In recent years the real Milnor fibrations has been extended in several directions in the settings of isolated and non-isolated singularities. Some results can be found for instance in \cite{AT, TYA, RS} and references.

\end{remark}

Milnor also provided information about the topology of the fiber of the fibration $(1),$ see Proposition \ref{PM1}, section 2. In fact, it follows from Lemma 11.4, page 100, that if the link
$K_{\epsilon}$ is not empty, then the fiber is $(p-2)-$connected. That is, for $p=2$ it is connected and for $p\geq 3$ it is simple connected.

\vspace{0.2cm}

More recently in \cite{ADN}, using tools from singularity theory, Morse theory and differential geometry,
the authors proved an extended Khinshiasvilli's formula \cite{Kh} for the Euler-Poincaré characteristic
of an isolated singularity of a polynomial map germ
$\displaystyle{f:(\mathbb{R}^{n},0)\to (\mathbb{R}^{p},0)}$,
$f(x)=(f_{1}(x),\dots , f_{p}(x))$, $\displaystyle{n\geq p \geq 2}$. Denote by $F_{f}$ the fiber of the
associated Milnor fibration $(1).$ The main result in \cite{ADN} is:

\begin{theorem}[\cite{ADN}]\label{T1}

Given a polynomial map germ $f$ as above, the following holds:

\begin{itemize}

\item[a)] If $n$ is even, then $\chi(F_{f})=1-deg_{0}(\nabla f_{1})$, where $deg_{0}(\nabla f_{1})$ is the topological degree of the mapping $\displaystyle{\frac{\nabla f_{1}}{\|\nabla f_{1}\|}:S_{\epsilon}^{n-1}\to S_{1}^{n-1}}$. Moreover, $deg_{0}(\nabla f_{1})=deg_{0}(\nabla f_{2})=\cdots = deg_{0}(\nabla f_{p})$.

\vspace{0.4cm}

\item[b)] If $n$ is odd, then $\chi(F_{f})=1.$ Moreover, $deg_{0}(\nabla f_{i})=0$, for all $i=1,\cdots, n.$

\end{itemize}

\end{theorem}

\vspace{0.2cm}

In \cite{Mi}, page 100, Milnor posed the following
question:

\vspace{0.2cm}

`` For which dimensions $n\geq p \geq 2$ do non-trivial example
exist ? "

\vspace{0.2cm}

Indeed Milnor did not specify what a `` trivial map '' should mean, but in a certain sense he proposed to call an isolated singularity polynomial mapping $f$ `` trivial '' if the fiber $F_{f}$ of its associated fiber bundle $(1)$ is diffeomorphic to the $\displaystyle{( n-p )-}$dimensional closed disk.

\vspace{0.2cm}

Using this concept as a definition, Church and Lamotke in \cite{CL} answered the above question in the following way.

\begin{theorem}[\cite{CL}, pags. 149--150]\label{T2}

\begin{itemize}
\item [a)] For $0\leq n-p\leq 2$, non-trivial examples occur precisely for the dimension $(n,p)=$
$(2,2)$, $(4,3)$ and $(4,2).$

\item [b)] For $n-p\geq 4$ non-trivial examples occur for all $(n,p).$
\item [c)] For $n-p=3$ non-trivial examples occur for $(5,2)$ and $(8,5).$ Moreover,
if the $(3-dimensional)$ Poincaré Conjecture is false, then there
are non-trivial examples for all $(n,p).$ If the Poincaré Conjecture
is true, all examples are trivial except $(5,2)$, $(8,5)$ and
possibly $(6,3).$

\end{itemize}
\end{theorem}

Furthermore, in \cite{CL} page 151, the authors proved an alternative
characterization of trivial map-germs (see below) using the branch set
$B_{f}$ of a map-germ $f$, i.e., the set where the map-germ fails to
be locally topologically equivalent to a projection.

\begin{proposition}[\cite{CL}, page 151]  For $n-p \neq 4,5$ $f$ is trivial if, and only if, $0\notin B_{f}.$

\end{proposition}

Also in \cite{CL} the authors found the pairs of dimensions $(n,p)$ for which $0$ is isolated in $f^{-1}(0)$, i.e., locally the level set $\{f^{-1}(0)\}=\{0\}.$ This is equivalent to the link $K_{\epsilon}=f^{-1}(0)\cap S_{\epsilon}^{n-1}$ being empty.

\begin{lemma}[\cite{CL}, Lemma 1, page 151]\label{L1}
If $0$ is an isolated point of $f^{-1}(0),$ then $n=p$, or $(n,p)=(4,3), (8,5), (16,9).$ If $n=p$ the singularity is trivial unless $k=2.$ The other cases are never trivial.

\end{lemma}

\vspace{0.2cm}

The aim of this paper is to use tools from algebraic topology, singularity theory and a formula for the Euler-Poincaré characteristic of the real Milnor fiber proved in \cite{ADN} (Theorem \ref{T1} above) to find conditions under which provide a (new) characterization of trivial map-germs. The main Theorem cover all
cases when the dimension of the Milnor fiber is $n-p=3$, i.e., the
item $(c)$ of Theorem \ref{T2} above. Our main result is:

\begin{theorem}\label{teorema 1}

Let $\displaystyle{f:(\mathbb{R}^{n},0)\to (\mathbb{R}^{p},0)}$,
$f(x)=(f_{1}(x),\dots , f_{p}(x))$ be a
polynomial map-germ with an isolated singularity, and suppose that $n-p=3$. Then, the following holds:

\begin{itemize}

\item[1)] If the pair $(n,p)=(6,3)$, then $f$ is trivial if and only if $\deg_{0} (\nabla f_{1} )=0$;

\item[2)] If the pair $(n,p)=(8,5)$, then $f$ is trivial if and only if the link $K_{\epsilon}$ is not empty;

\item[3)] If the pair $(n,p)=(5,2)$, then $f$ is trivial if and only if $\pi_{1}(F_{f})=0$
 (i.e. the respective Milnor fiber  $F_{f}$ is simply connected).

\end{itemize}
\end{theorem}

It is worth pointing out that in \cite{ADN}, section $4$, the authors proved some formulae in order to describe geometrically/topologically the real Milnor fiber and a characterization of trivial map-germ for the pair $(4,2)$. For this they used a condition analogous to that in the Theorem \ref{teorema 1} $(1)$, (see \cite{ADN}, Corollary 4.5). They proved that an isolated singularity polynomial map-germ $\displaystyle{f=(f_{1},f_{2}):(\mathbb{R}^{4},0)\to (\mathbb{R}^{2},0)}$ is trivial if and only if $\deg_{0}(\nabla f_{1})=0.$

\section{Set up and preliminary results}

In this section we state and prove some results from algebraic
topology in order to prove our main result. We focus on results related to the topology of three
dimensional compact manifolds.

Next result is a very well known
fact from algebraic topology, and we include it here for the sake of completeness. The prove follows from the additive property of Euler-Poincaré characteristic, (see also
\cite{Vick}, Theorem 6.38, p. 180).

\begin{theorem}[]\label{euler caracteristic}
If $M$ is a compact to\-po\-lo\-gi\-cal $n$-ma\-ni\-fold with
boundary $\partial M$, then the Euler characteristic of $\partial M$ is even. More precisely,
\begin{eqnarray*}\chi(\partial M)=\left\{\begin{array}{ll}
2\chi(M) & \text{if $n$ is odd,}\\
0& \text{if $n$ is even.}
\end{array}\right.
\end{eqnarray*}

\end{theorem}

\begin{proposition}\label{proposition1} If $M$ is a connected, compact and orientable
$n$-manifold with non empty boundary $\partial M$, then
$H_i(M;\mathbb{Z})=0$, for $i\geq n$.
\end{proposition}

\begin{proof} Since $\stackrel{\circ}{M}$ is a connected non-compact
$n$-manifold without boun\-dary and the inclusion map
$i:\stackrel{\circ}{M}\hookrightarrow M$ is an homotopy equivalence,
by \cite[Proposition 3.29, p. 239]{Hatcher}, we have
\begin{eqnarray}\nonumber H_{i}(M;\mathbb{Z})\cong
H_{i}(\stackrel{\circ}{M};\mathbb{Z})=0,\hspace{0.2cm}{\rm
for}\hspace{0.2cm} i\geq n.
\end{eqnarray}
\end{proof}

\begin{lemma}\label{lema 1} Let $M$ be a compact, connected and orientable $3$-manifold with boundary $\partial M$.
If $H_{1}(M,\mathbb{Z})=0$, then $\partial M$ is a disjoint union of $2$-spheres.
\end{lemma}

\begin{proof} Consider the following piece of the long exact
sequence in cohomology for the pair $(M,\partial M)$:

\begin{eqnarray}\xymatrix{\cdots\rightarrow H^1(M,\mathbb{Z})\rightarrow
H^1(\partial M; \mathbb{Z})\rightarrow H^2(M,\partial
M;\mathbb{Z})\rightarrow\cdots}\label{sequencia 1}\end{eqnarray}

Since $M$ is compact and orientable, by Poincaré-Lefschetz duality
theorem (see \cite{Vick}, Theorem 6.25, p. 171) we have
\begin{eqnarray}\nonumber H^{2}(M,\partial M;\mathbb{Z})\cong
H_{1}(M;\mathbb{Z})=0.
\end{eqnarray}
On the other hand, $H^{1}(M;\mathbb{Z})\cong Hom
(H_{1}(M;\mathbb{Z}), \mathbb{Z})=0$, therefore, from (\ref{sequencia 1}) it
follows that $H^{1}(\partial M;\mathbb{Z})=0$ and by using the
Poincaré duality for the orientable closed $2$-manifold $\partial
M$,
\begin{eqnarray}\nonumber 0=H^{1}(\partial M;\mathbb{Z})\cong H_{1}(\partial M;\mathbb{Z}).\end{eqnarray}

\noindent By the standard classification of the compact
$2$-manifolds, we conclude that $\partial M$ is a disjoint union of
$2$-spheres.
\end{proof}

\begin{lemma}\label{lema 2} Let $M$ be a compact, connected and orientable $3$-manifold with boundary $\partial M$. Then the
following statements are equivalent:

\vspace{0.2cm}

{\rm(i)} $M$ is contractible.

\vspace{0.2cm}

{\rm(ii)} $M$ is simply-connected and has the same homology of a
point.

\vspace{0.2cm}

{\rm(iii)} $M$ is simply-connected and $\partial M$ is a 2-sphere.

\vspace{0.2cm}

{\rm(iv)}  $M$ is diffeomorphic to a 3-disk.
\end{lemma}

\begin{proof} It is follow from \cite[Corollary 10.11, p. 479]{Bredon} that $(i)$ is equivalent to $(ii).$

\noindent For (ii)$\Leftrightarrow $(iii), since $M$ is simply-connected,
from Lemma \ref{lema 1}, $\partial M$ is a disjoint union of
$2$-spheres. Therefore, if $M$ has the same homology of a point, we
have $\chi(M)=1$, and by Theorem \ref{euler caracteristic},
\begin{eqnarray}\nonumber \chi(\partial M)=2\chi(M)=2,
\end{eqnarray}
that is, $\partial M$ is a $2$-sphere.

Conversely, since $M$ is a simply-connected space (in particular, $H_{1}(M;\mathbb{Z})=0$) and from Proposition
\ref{proposition1}, $H_{j}(M;\mathbb{Z})=0$ for $j\geq 3$, it is
sufficient to show that  $H_{2}(M;\mathbb{Z})=0$. Indeed, we have
$H^{1}(M;\mathbb{Z})\cong Hom (H_{1}(M;\mathbb{Z}), \mathbb{Z})=0$
and by Poincaré-Lefschetz duality theorem
\begin{eqnarray*}0\cong H^{1}(M;\mathbb{Z})\cong H_{2}(M,\partial M;\mathbb{Z}).
\end{eqnarray*}

Now, by  \cite[Theorem 6.24, p. 169]{Vick},
there is a unique fundamental class $z\in H_{3}(M,\partial M;\mathbb{Z})$
such that $\Delta(z)$ is a fundamental class of $\partial M$,
where $\Delta:H_{3}(M,\partial M;\mathbb{Z}) \to H_{2}(\partial M;\mathbb{Z})$
is the connecting homomorphism.
By assumptions $\partial M$ is a 2-sphere,
 we conclude that $\Delta$ is onto and considering a
particular  part of the long exact sequence in homology for the pair $(M,\partial M)$
\begin{eqnarray}\nonumber\label{sequencia 2}\xymatrix{
H_{3}(M,\partial M)\stackrel{\Delta}{\twoheadrightarrow}
H_{2}(\partial M)\rightarrow H_{2}(M)\stackrel{}{\rightarrow}
H_{2}(M,\partial M)=0, }\\\nonumber\end{eqnarray}
it follows that $H_{2}(M;\mathbb{Z})=0$.

\vspace{0.2cm}

For (i)$\Leftrightarrow$ (iv),  let $N=M\cup_{\partial M} M$ be the
``double" of $M$ and let $M_{1}$ and $M_{2}$ be two copies of $M$ in
$N$. We have that $N$ is a closed orientable $3$-manifold. By
\cite[Proposition 3.42, p. 253]{Hatcher} we can consider $U$ the
union of $M_{1}$ and a collar neighborhood of $M_{2}$; similarly we
can consider $V$ the union of $M_{2}$ and a collar neighborhood of
$M_{1}$. Since the pathwise connected spaces $U,V$ and $U\cap V$ are
homotopy equivalents to $M_{1}$, $M_{2}$ and $\partial M$,
respectively and $\partial M$ is simply connected, by Van Kampen
Theorem (see \cite[Corollary 4.27, p. 116]{Vick}), it follows that
\begin{eqnarray}\nonumber \pi_{1}(N)\cong \pi_{1}(U)*
\pi_{1}(V)\cong \pi_{1}(M_{1})* \pi_{1}(M_{2}).
\end{eqnarray}

Therefore, since by assumptions  $M_{1}=M_{2}=M$ are contractible, we have
$\pi_{1}(N)=0$. Thus, by Poincaré Theorem we conclude that $N$ is
diffeomorphic to a $3$-sphere and consequently, $M$ is diffeomorphic
to a $3$-disk. The converse it is obvious.
\end{proof}

\vspace{0.3cm}

Next proposition follow from \cite{Mi}, Lemma 11.4, page 100, which gives an important information about the connectedness of the real Milnor fibre.

\begin{proposition}\label{PM1} Let $\displaystyle{f:(\mathbb{R}^{n},0)\to (\mathbb{R}^{p},0)}$,
$f(x)=(f_{1}(x),\dots , f_{p}(x))$ be an isolated singularity
polynomial map-germ, with $n\geq p\geq 2$. If the link $K_{\epsilon}$ is not empty, then the fiber $F_{f}$ of the fibration $(1)$ is $(p-2)-$connected.

\end{proposition}

\begin{proof} See the first two paragraphs on the proof of Lemma 11.4, pages 100--101, in \cite{Mi}.  \end{proof}


\section{Proof of the main result}

In this section we are assuming that $\displaystyle{f:(\mathbb{R}^{n},0)\to (\mathbb{R}^{p},0)}$,
$f(x)=(f_{1}(x),\dots , f_{p}(x))$ is an isolated singularity
polynomial map-germ, with $n-p=3$.
For simplicity, we split the prove of Theorem \ref{teorema 1} into three parts.
Below we start with the proof of the item (1).
Recall that, in the case $(6,3)$ it follows from Lemma. \ref{L1} that the link $K_{\epsilon},$ which is diffeomorphic to the boundary of the Milnor fiber $\partial F_{f},$ is not empty.

\subsection{ \underline{The (6,3) case}}

Let $f: (\mathbb{R}^{6},0)\to (\mathbb{R}^{3},0)$ be with isolated
singu\-la\-ri\-ty at origin.

\begin{proof} ({\bf Proof of Theorem \ref{teorema 1}, item (1)}) : If $f$ is trivial, $F_{f}$ is diffeomorphic to a $3$-disk, so the Euler Characteristic of $F_{f}$ satisfies $\chi(F_{f})=1$. As $\chi(F_{f})=1-\deg_{0}\nabla f_{1}$ (Theorem 1.2 $(a)$) we conclude that $\deg_{0}\nabla f_{1}=0$.

Conversely, since $p=3$ the Milnor fiber $F_{f}$ is simply-connected, so in particular
$H_{1}(F;\mathbb{Z})=0$. Now, from Lemma \ref{lema 1} we have that
$\partial F_{f} $, which is diffeomorphic to $K_{\epsilon}$, is a disjoint union of $2-$spheres. Therefore, if $deg_{0}\nabla f=0$, then $\chi(F_{f})=1$, and
\begin{eqnarray}\nonumber \chi(K_{\epsilon})=\chi(\partial F_{f})=2\chi(F_{f})=2.
\end{eqnarray}
Hence, the link $K_{\epsilon}$ is diffeomorphic to one $2$-sphere and by Lemma \ref{lema 2}, we
have that $F_{f}$ must be diffeomorphic to a $3$-disk.
\end{proof}

\begin{proposition} For the pair $(n,p)=(6,3)$, the following
are equi\-va\-lent:

\vspace{0.1cm}

{\rm(1)} $f$ is trivial;

\vspace{0.1cm}

{\rm(2)} The link $K_{\epsilon}$ is connected.

\end{proposition}

\begin{proof} The implication $(1)\Rightarrow (2)$ is trivial.
In fact, if we assume that $f$ is trivial, the boundary
of the fiber $\partial F_{f}$ is diffeomorphic to the $2-$sphere, and so is
connected.

For the converse, if we assume that the link is connected, we know
that $K_{\epsilon}$ must be diffeomorphic to one single copy of a
$2$-sphere and again by Lemma \ref{lema 2} $(iii)$, $F$
is diffeomorphic to a $3$-disk. \end{proof}

We can know summarize all these results in the following way.

\begin{corollary}\label{corolario 3.2}

Given $f: (\mathbb{R}^{6},0)\to (\mathbb{R}^{3},0),$
$f(x)=(f_{1}(x),f_{2}(x),f_{3}(x)),$ an isolated singularity
polynomial map-germ, the following statements are equivalent:

\vspace{0.2cm}

{\rm(1)} $f$ is trivial;

\vspace{0.1cm}

{\rm(2)} The link $K_{\epsilon}$ is connected;

\vspace{0.1cm}

{\rm(3)} $\deg_{0}(\nabla f_{1})=0.$

\end{corollary}

\vspace{0.2cm}

\begin{remark}\rm The statement (2) of Corollary \ref{corolario 3.2} does
not cha\-racte\-rizes in general the trivial fibration. For example, in the
case of the pair $(4,2)$, if we consider $f(x,y)=x^{2}+ y^{3},$ then the link is the $(2,3)-$torus knot and so it is connected, and the open fiber is diffeomorphic to the torus without an open disc removed.

\end{remark}

\begin{proposition}

Given $f: (\mathbb{R}^{6},0)\to (\mathbb{R}^{3},0),$
an isolated singularity polynomial map-germ,
then the Euler-Poincaré characteristic of the fiber $\chi(F_{f})$ is precisely the number of spheres on its boundary $\partial F_{f}$.

\end{proposition}

\begin{proof}
\noindent Since $F_{f}$ is a $3-$dimensional, simply
connected and compact manifold with boundary given by $d-$disjoint
copies of $2-$dimensional spheres, then each sphere on $\partial
F_{f}$ contributes with $2$ to its Euler cha\-racteristic. Hence,
$$2.d=\displaystyle{\chi(\partial F_{f})=2.\chi(F_{f})}.$$

Therefore, the Euler characteristic of the fiber $F_{f}$ is the number of spheres in the boundary.

\end{proof}

\subsection{\underline{The (8,5) case}}

Consider $f: (\mathbb{R}^{8},0)\to (\mathbb{R}^{5},0)$ a polynomial map-germ with an isolated
singularity at origin. From Lemma \ref{L1} we
know that for the pair $(8,5)$ the link may be empty. Of course that,
for the purpose of characterization of trivial mapping, a necessary
condition is that the link be not empty.

\begin{proof} ({\bf Proof of Theorem \ref{teorema 1}, item (2)}):
The first implication is obvious. Conversely, since the link
$K_{\epsilon}\not=\emptyset$ and $p=5$, the Milnor
fiber $F_{f}$ is a $3$-connected (and so is simple connected), compact
$3$-manifold with boundary. From Hurewicz Theorem, $F_{f}$ must have
the homology type of a point, and so by Lemma \ref{lema 2}
$(ii),$ $F_{f}$ is diffeomorphic to a $3$-disk.
\end{proof}

\begin{remark}
By a similar argument to that used in the proof of Theorem \ref{teorema 1} $(1),$
it follows that the condition $\deg_{0}(\nabla f_{1})=0$ also characterizes the
triviality of $f$ in this case. Hence, we have:
\end{remark}
\begin{corollary} For the pair $(n,p)=(8,5)$, the following statements
are equivalent:

\vspace{0.2cm}

{\rm(1)} $f$ is trivial;

\vspace{0.1cm}

{\rm(2)} The link $K_{\epsilon}\not=\emptyset$;

\vspace{0.1cm}

{\rm(3)} $\deg_{0}(\nabla f_{1})=0$.

\end{corollary}

\subsection{\underline{The case (5,2)}} Now consider $f: (\mathbb{R}^{5},0)\to (\mathbb{R}^{2},0)$
be a polynomial map-germ with an isolated singularity at origin. In this case, we only have
that the Milnor fiber $F_{f}$ is connected and the link is
not empty.

\begin{proof} ({\bf Proof of Theorem \ref{teorema 1}, item (3)}) :
Again, for the first implication there is nothing to prove, since the closed disk is simple connected. For the converse, if we assume that $F_{f}$ is simply-connected, i.e., $\pi_{1}(F_{f})=0,$ then we know from Lemma \ref{lema 1} that
$\partial F_{f} \approx K_{\epsilon} $ is diffeomorphic to a disjoint union of finite many 2-spheres. In another hand, since the source space have odd dimension, then by
 Theorem \ref{T1} (b), see \cite{ADN}, we have $\chi(F_{f})=1$ and so
\begin{eqnarray}\nonumber \chi(K_{\epsilon})=\chi(\partial F_{f})=2.
\end{eqnarray}
Therefore, the link $K_{\epsilon}$ must be a $2$-sphere and by Lemma \ref{lema 2}, we
have that $F_{f}$ is diffeomorphic to a $3$-disk.
\end{proof}

\begin{remark} In \cite{Lo}, page 421, E. Looijenga remarked that we can
use the example $f:(\mathbb{C}^{2},0)\to (\mathbb{C},0)$,
$f(x,y)=x^{2}+y^{3}$, which is not trivial (see for
instance \cite{ADN}), as a real polynomial map-germ
from $(\mathbb{R}^{4},0) \to (\mathbb{R}^{2},0)$ to guarantee the existence of
 non-trivial polynomial map-germ from $(\mathbb{R}^{5},0)
\to (\mathbb{R}^{2},0).$ Therefore, as a by-product of our
characterization above, we can conclude that its respective Milnor
fiber can not be simply connected.
\end{remark}

\begin{remark} In the case $(5,2)$, if we suppose that
$H_{1}(F_{f};\mathbb{Z})=0$, from Lemma \ref{lema 1} and Theorem
\ref{T1} (b), then the link $K_{\epsilon}$ is a
$2$-sphere. However, we cannot ensure that $F_{f}$ is diffeomorphic
to a $3$-disk, as the following example shows: by \cite[Theorem 8.10
(Poincaré), page 353]{Bredon}, there is a compact 3-manifold $W$
having the homology groups of the sphere $S^{3}$ but which is not
simply-connected. So, $W\setminus \stackrel{\circ}{D^{3}}$ is a
compact 3-manifold whose boundary is a 2-sphere,
which is not simply-connected.
\end{remark}

\end {document}